\documentclass[12pt]{amsart}
\usepackage{amsmath}
\usepackage{amsxtra}
\usepackage{amscd}
\usepackage{amsthm}
\usepackage{amsfonts}
\usepackage{amssymb}
\usepackage{eucal}
\usepackage{mathrsfs}
\usepackage{hyperref}
\usepackage{pstricks,pst-plot}

\headsep=1cm \numberwithin{equation}{section}
\theoremstyle{definition}
\newtheorem{theorem}{Theorem}[section]

\newtheorem{proposition}[theorem]{Proposition}
\newtheorem{corollary}[theorem]{Corollary}
\newtheorem{remark}[theorem]{Remark}

\newtheorem{convention and reminder}[theorem]{Convention and Reminder}
\newtheorem{convention and remark}[theorem]{Convention and Remark}
\newtheorem{definition and remark}[theorem]{Definition and Remark}

\newtheorem{reminders and definition}[theorem]{Reminders and Definition}

\newtheorem{notation and remarks}[theorem]{Notation and Remarks}
\newtheorem{notation and remark}[theorem]{Notation and Remark}
\newtheorem{example}[theorem]{Example}

\def\N{{\mathbb N}}

\def\P{{\mathbb P}}

\newcommand\Ker{\operatorname{\Ker}}

\begin{document}

\title{On hypersurfaces containing projective varieties}

\author{Euisung Park}


\begin{abstract}
Classical Castelnuovo's Lemma shows that the number of linearly independent quadratic equations of a nondegenerate irreducible projective variety of codimension $c$ is at most ${{c+1} \choose {2}}$ and the equality is attained if and only if the variety is of minimal degree. Also a generalization of Castelnuovo's Lemma by G. Fano implies that the next case occurs if and only if the variety is a del Pezzo variety. For curve case, these results are extended to equations of arbitrary degree respectively by J. Harris and S. L'vovsky. This paper is intended to extend these results to arbitrary dimensional varieties and to the next cases.
\end{abstract}

\maketitle
\thispagestyle{empty}

\section{Introduction}
\noindent Let $X \subset \P^r$ be a nondegenerate irreducible projective variety defined over an algebraically closed field $K$ of arbitrary characteristic. A basic invariant of $X$ is the number of linearly independent hypersurfaces of degree $m$ containing $X$ for each $m \geq 2$. Throughout this paper, we will denote this number by $a_m (X)$. That is, $a_m (X) = h^0 (\P^r , \mathcal{I}_X (m))$ where $\mathcal{I}_X$ is the sheaf of ideals of $X$ in $\P^r$. The aim of this paper is to find an upper bound for $a_m (X)$ and to investigate the borderline cases.

As S. L'vovsky indicates in \cite{L}, the first results in this direction are due to G. Castelnuovo and G. Fano who proved respectively that a finite set $\Gamma \subset \P^c$ of $d$ points in linearly general position should lie on
\begin{enumerate}
\item[(i)] (Classical Castelnuovo's Lemma) a rational normal curve if $d \geq 2c+3$ and $h^0 (\P^c , \mathcal{I}_{\Gamma} (2))={{c} \choose {2}}$, and
\item[(ii)] a linearly normal curve of arithmetic genus one if $d \geq 2c+5$ and $h^0 (\P^c , \mathcal{I}_{\Gamma} (2))={{c} \choose {2}}-1$.
\end{enumerate}
Nowadays (ii) was rediscovered by D. Eisenbud and J. Harris \cite{H2}. These results imply the following

\begin{theorem}[Theorem 1.2 and Theorem 1.6 in \cite{L}]\label{thm:quadratic}
Let $X \subset \P^r$ be a nondegenerate irreducible projective variety of codimension $c \geq 2$ and degree $d$. Then

\smallskip

\renewcommand{\descriptionlabel}[1]%
             {\hspace{\labelsep}\textrm{#1}}
\begin{description}
\setlength{\labelwidth}{13mm}
\setlength{\labelsep}{1.5mm}
\setlength{\itemindent}{0mm}

\item[(1)] (G. Castelnuovo) $X$ is contained in at most ${{c+1} \choose {2}}$ linearly independent quadrics; the equality is attained if and only if $d=c+1$.
\item[(2)] (G. Fano) $X$ is contained in exactly ${{c+1} \choose {2}}-1$ linearly independent quadrics if and only if its one-dimensional general linear sections are linearly normal curves of arithmetic genus one.
\end{description}
\end{theorem}

\noindent Theorem \ref{thm:quadratic} was  reproved by F. L. Zak \cite{Z} whose proofs make extensive use of secant varieties. For projective curves, the following is known:

\begin{theorem}[Theorem 1.4 and Theorem 1.7 in \cite{L}]\label{thm:curves}
Let $C \subset \P^r$, $r \geq 3$, be a nondegenerate projective integral curve and let $m \geq 2$ be an integer. Then

\smallskip

\renewcommand{\descriptionlabel}[1]%
             {\hspace{\labelsep}\textrm{#1}}
\begin{description}
\setlength{\labelwidth}{13mm}
\setlength{\labelsep}{1.5mm}
\setlength{\itemindent}{0mm}

\item[(1)] (J. Harris, \cite{H1}) $a_m (C) \leq {{r+m} \choose {r}}- (mr+1)$.
Also the equality is attained if and only if $C$ is a rational normal curve.
\item[(2)] (S. L'vovsky) If $C$ is not a rational normal curve, then
\begin{equation*}
a_m (C) \leq {{r+m} \choose {r}}- m(r+1).
\end{equation*}
Also the equality is attained if and only if $C$ is a linearly normal curve of arithmetic genus one.
\end{description}
\end{theorem}

\noindent We refer the reader to \cite{L} for a nice survey about the problem mentioned at the beginning. In this paper we will extend Theorem \ref{thm:quadratic} and Theorem \ref{thm:curves} to arbitrary $n$ and $m$.\\

To state our main result precisely, we require some notation and remarks. For integers $n \geq 1$ and $c \geq 2$, let \begin{equation*}
\Xi_{n,c} := \{ X \subset \P^{n+c}~|~ X:~\mbox{$n$-dimensional nondegenerate projective variety} \}
\end{equation*}
be the set of all nondegenerate irreducible projective varieties of dimension $n$ in $\P^{n+c}$. For each $m \geq 2$, we regard $a_m$ as the function from $\Xi_{n,c}$ to the set $\N_0$ of nonnegative integers. Obviously the image
\begin{equation*}
A_{n,c,m} := \{ a_m (X) ~|~ X \in \Xi_{n,c} \}
\end{equation*}
of $a_m$ is a finite subset of $\N_0$. Thus we can define
\begin{equation*}
\delta_{n,c,m} (k) := \mbox{the $k$th largest member of $A_{n,c,m}$}
\end{equation*}
for $1 \leq k \leq |A_{n,c,m}|$. Keeping these notations in mind, we can reformulate the problem outlined at the first paragraph of this section in the following form:\\
\begin{enumerate}
\item[] {\bf Problem A.} For each $k \geq 1$, determine the value of $\delta_{n,c,m} (k)$ and find all $X \in \Xi_{n,c}$ satisfying $a_m (X) = \delta_{n,c,m} (k)$. \\
\end{enumerate}

\noindent Theorem \ref{thm:quadratic} and Theorem \ref{thm:curves} provides an answer for $m=2$ and $n=1$, respectively. For example, they mean that
$$ \begin{cases} \delta_{n,c,2} (1) & \quad = \quad {{c+1} \choose {2}},  \\
\delta_{n,c,2} (2)                  & \quad = \quad {{c+1} \choose {2}}-1,  \\
\delta_{1,c,m} (1)                  & \quad = \quad {{c+1+m} \choose {m}}- \{m(c+1)+1\} \quad \mbox{and}\\
\delta_{1,c,m} (2)                  & \quad = \quad {{c+1+m} \choose {m}}- m(c+2). \end{cases}$$
Now, we define the following functions on the positive integers $n$, $c$ and $m$ where $1 \leq t \leq n+1$:

\begin{equation*}
F(n,c,m) = {{m+n+c} \choose {n+c}} - \{ (c+1) {{m+n-1} \choose {n}} + {{m+n-1} \choose {n-1}} \}
\end{equation*}

\begin{equation*}
G_t (n,c,m) = {{m+n+c} \choose {n+c}} \quad \quad \quad \quad \quad \quad \quad \quad \quad \quad \quad \quad \quad \quad \quad \quad \quad \quad
\end{equation*}
\begin{equation*}
\quad \quad \quad \quad \quad - \{ (c+2) {{m+n-1} \choose {n}} + {{m+n-1} \choose {n-1}} - {{m+t -3} \choose {t-2}} \}
\end{equation*}
\smallskip

\noindent Here ${{a} \choose {b}}=0$ if $b < 0$ or $a<b$. These integer-valued functions come from some projective varieties of low degree. Namely, for an $n$-dimensional projective irreducible variety $X \subset \P^{n+c}$ of degree $d$ it holds that
$$a_m (X) = \begin{cases} F(n,c,m) & \mbox{if $d=c+1$, and} \\
                          G_t (n,c,m) & \mbox{if $d=c+2$ and ${\rm depth}(X)=t$.}  \end{cases}$$
We denote by ${\rm depth}(X)$ the arithmetic depth of the homogeneous coordinate ring of $X$. The first case is probably well-known. For the second case, we refer to Theorem A and Theorem B in \cite{HSV} (see also Theorem 2.2 in \cite{N}).

Our main result is the following

\begin{theorem}\label{thm:main}
Let $X \subset \P^{n+c}$, $c \geq 2$, be an $n$-dimensional nondegenerate projective variety of degree $d$ and let $m \geq 2$ be an integer.

\smallskip

\renewcommand{\descriptionlabel}[1]%
             {\hspace{\labelsep}\textrm{#1}}
\begin{description}
\setlength{\labelwidth}{13mm}
\setlength{\labelsep}{1.5mm}
\setlength{\itemindent}{0mm}

\item[(1)] $\delta_{n,c,m} (1)=F(n,c,m)$. Also $a_m (X)= F (n,c,m)$ if and only if $d=c+1$.
\item[(2)] $\delta_{n,c,m} (2)=G_{n+1} (n,c,m)$. Also $a_m (X)= G_{n+1} (n,c,m)$ if and only if $d=c+2$ and ${\rm depth}(X)=n+1$.
\item[(3)] $\delta_{n,c,m} (3)=G_n (n,c,m)$. Also for $m \geq 3$, $a_m (X) = G_n (n,c,m)$ if and only if $d=c+2$ and ${\rm depth}(X)= n$.
\item[(4)] For $c \geq 3$, $a_2 (X)= G_n (n,c,2)$ if and only if either
$d=c+2$ and ${\rm depth}(X)= n$ or else $d=c+3$ and ${\rm depth}(X)=n+1$.
\end{description}
\end{theorem}
\smallskip

For the proof of this result see Theorem \ref{thm:firstandsecond}. Also see Remark \ref{rem:4.7} in which we discuss about $\delta_{n,c,m} (4)$ for $m \geq 3$. Theorem \ref{thm:main} says that if $m \geq 3$, then some positive integers $\leq F(n,c,m)$ are not contained in $A_{n,c,m}$ since
\begin{equation*}
\delta_{n,c,m} (1) - \delta_{n,c,m}(2) = {{m+n-2} \choose {n}} > 1.
\end{equation*}
Note that Theorem \ref{thm:main}.(1) and (2) imply Theorem \ref{thm:quadratic} and Theorem \ref{thm:curves} since
\begin{equation*}
F(n,c,2) = {{c+1} \choose {2}} \quad \mbox{and} \quad  G(n,c,2) = {{c+1} \choose {2}}-1
\end{equation*}
are respectively Castelnuovo's and Fano's bounds and
\begin{equation*}
F(1,c,m) = {{r+m} \choose {r}}- (mr+1) \quad \mbox{and} \quad  G(1,c,m) = {{r+m} \choose {r}}- m(r+1)
\end{equation*}
are respectively Harris's and L'vovsky's bounds where $r=c+1$. Also Theorem \ref{thm:main}.(3) and (4) extend Theorem \ref{thm:quadratic} and Theorem \ref{thm:curves} to the next cases, respectively. More precisely, we have the following corollaries:

\begin{corollary}\label{cor:quadratic}
Let $X \subset \P^{n+c}$ be an $n$-dimensional nondegenerate projective variety of degree $d$. If $c \geq 3$, then
\begin{equation*}
a_2 (X) = {{c+1} \choose {2}}-2
\end{equation*}
if and only if either $d=c+2$ and ${\rm depth}(X)=n$ or else $d=c+3$ and ${\rm depth}(X)=n+1$.
\end{corollary}

\begin{corollary}\label{cor:curves}
Let $C \subset \P^r$ be a nondegenerate projective integral curve and let $m \geq 2$ be an integer. If $r = 3$ and $m \geq 3$ or if $r \geq 4$, then
\begin{equation*}
a_m (C) = {{r+m} \choose {r}}- m(r+1)-1
\end{equation*}
if and only if either $C$ is a smooth rational curve of degree $r+1$
or else $m=2$ and $C$ is a linearly normal curve of arithmetic genus two.
\end{corollary}

Briefly speaking, the proof of Theorem \ref{thm:main} is based on the hyperplane section method and the induction argument on $n$ and $m$. To be precise, let $\Gamma \subset \P^c$ and $C \subset \P^{c+1}$ be respectively general zero-dimensional and one-dimensional linear sections of $X$. Since $\Gamma$ spans $\P^c$, we may assume that it contains the set $\Gamma_0$ of the $(c+1)$ coordinate points of $\P^c$. Then it follows that
\begin{equation}\label{eq:1.1}
a_2 (X) \leq h^0 (\P^c , \mathcal{I}_{\Gamma} (2)) \leq h^0 (\P^c , \mathcal{I}_{\Gamma_0} (2)) = {{c+1} \choose {2}}.
\end{equation}
Note that the upper bound ${{c+1} \choose {2}}$ of the number of linearly independent quadratic equations of $X$ is obtained by using the fact that $\Gamma_0$ is $3$-regular. In Proposition \ref{prop:semiuniform}, we generalize this elementary result by showing that $\Gamma$ contains a subset $\Gamma'$ of ${\rm min} \{d,2c+1\}$ points which is $3$-regular and spans $\P^c$. This fact is well-known if ${\rm char}(K)=0$ but its proof demands more effort if ${\rm char}(K)$ is positive. From (1.1) one can naturally pose the following\\
\begin{enumerate}
\item[] {\bf Problem B.} For each $n \geq 1$, $c \geq 2$ and $k \geq 1$, classify all $n$-dimensional nondegenerate projective varieties $X \subset \P^{n+c}$ with
\begin{equation}\label{eq:1.2}
a_2 (X) ={{c+1} \choose {2}}+1-k.
\end{equation}
\end{enumerate}
\smallskip

\noindent Concerned with this problem, Proposition \ref{prop:semiuniform} enables us to show that if $k \leq c$, then (1.2) implies $d \leq c+k$ (see Corollary \ref{cor:quadraticinequality}). Therefore Problem B is closely related to the classification of varieties whose degree is at most $c+k$. Along this line, we study in Section 3 the deficiency module of projective integral curves $C \subset \P^{c+1}$ whose degree is at most $2c$. The most interesting result throughout this section is that when $C$ is not linearly normal, the integers
\begin{equation*}
h^1 (\P^{c+1},\mathcal{I}_C (j)),\quad 1 \leq j \leq {\rm reg}(C)-1,
\end{equation*}
form a strictly decreasing sequence (see Theorem \ref{thm:monotonousfunction}). From this result, we obtain a satisfactory answer for Problem B in curve case (see Theorem \ref{thm:curvehomologicaltype} and Theorem \ref{thm:quadric-curve}).

Section 4 is devoted to give a proof of Theorem \ref{thm:main} by combining the results in the previous two sections and some known facts about varieties of low degree.

In Section 5, we investigate some projective invariants of quadratic embedding defined by F. L. Zak\cite{Z} for varieties listed in Theorem \ref{thm:main}. It turns out that those varieties are characterized by their projective invariants of quadratic embedding (see Theorem \ref{thm:mainquadrics}).

Finally in Section 6, we solve Problem A completely when $n=1$ and $m \geq c$ (see Theorem \ref{thm:curvehigherdegree}).   Essentially this is possible because Theorem \ref{thm:monotonousfunction} guarantees a uniform cohomological behavior  of projective curves for all $m \geq c$.

\begin{remark}
Let $X \subset \P^{n+c}$ be as in Theorem \ref{thm:main}. \\
(1) If $d=c+1$, then $X$ is called a \textit{variety of minimal degree}. All such varieties have been completely classified by del Pezzo in the case of surfaces and by Bertini in the general case (cf. \cite{EH}). A variety of minimal degree is a cone over the Veronese surface in $\P^5$ or a rational normal scroll. It is arithmetically Cohen-Macaulay and cut out by quadrics.\\
(2) If $d=c+2$, then $X$ is called a \textit{variety of almost minimal degree} due to Brodmann-Schenzel\cite{BS3}. T. Fujita\cite{F2} has a satisfactory classification theory of those varieties. A variety $X$ of almost minimal degree is either normal and linearly normal or else obtained by projecting a variety of minimal degree. In the former case, ${\rm depth}(X)=n+1$ and $X$ is called a \textit{normal del Pezzo variety}. In the latter case, several basic properties (including the Hilbert function) of $X$ are firstly given in \cite{HSV}. Also the defining equations of $X$ and the syzygies among them are investigated in \cite{LP} and \cite{P}. Recently, a very detailed description of $X$ in terms of the projection map is obtained in \cite{BP} and \cite{BS3}. $X$ is cut out by quadrics if ${\rm depth}(X)=n+1$. But it may not be cut out by quadrics if ${\rm depth}(X) \leq n$ (e.g. Theorem 1.3 in \cite{P}).\\
(3) Varieties with $d=c+3$ are not yet completely classified. In this direction, we refer the reader to Section 10 in \cite{F2}, in which the classification of linearly normal smooth case is provided. It should be noted that if $c \geq 3$ and ${\rm depth}(X)=n+1$, then $X$ is cut out by quadrics.
\end{remark}

\noindent {\bf Acknowledgement.} This research was supported by Basic Science Research Program through the National Research Foundation of Korea(NRF) funded by the Ministry of Education, Science and Technology(2010-0007329). \\

\section{Finite sets in linear semi-uniform position}

\noindent Due to E. Ballico\cite{Ba}, a finite subset $\Gamma \subset \P^c$ is called in \textit{linear semi-uniform position} if it spans $\P^c$ and there are integers $\nu (i,\Gamma)$, $0 \leq i \leq c$, such that every $i$-plane $L$ in $\P^c$ spanned by linearly independent $i+1$ points of $\Gamma$ contains exactly $\nu (i,\Gamma)$ points of $\Gamma$. Thus $\Gamma$ is in \textit{linear general position} in the sense that any $c$ points in $\Gamma$ are linearly independent if and only if $\nu(c-1,\Gamma)=c$.

When $C \subset \P^{c+1}$ is a nondegenerate projective integral curve and $\Gamma \subset \P^c$ is its generic hyperplane section, $\Gamma$ is in linear semi-uniform position since every symmetric product of $C$ is irreducible. Furthermore, it is known that $\Gamma$ is in linear general position if ${\rm char}~K =0$ or if ${\rm chark}~K >0$, $c \geq 3$ and $C$ is smooth. Due to  E. Ballico\cite{Ba}, we say that $C$ is \textit{very strange} if $\Gamma$ fails to be in linear general position. There does exist very strange curves (cf. \cite[Example 1.2]{Ra})

A critical difference between the linear general position and the linear semi-uniform position is that if $\Gamma \subset \P^c$ is in linear general position then any subset $\Gamma' \subset \Gamma$ of more than $(c+1)$ points is still in linear general position while the linear semi-uniform position property may not be satisfied by a subset of $\Gamma$ if $\Gamma \subset \P^c$ is in linear semi-uniform position but fails to be in linear general position.

It is a well-known and elementary fact that if $\Gamma \subset \P^c$ is a finite set of $d (\geq 2c+1)$ points in linear general position, then any subset of $(2c+1)$ points of $\Gamma$ is $3$-regular. We generalize this fact to finite sets in linear semi-uniform position.

\begin{proposition}\label{prop:semiuniform}
Let $\Gamma \subset \P^c$ be a finite set of $d \geq 2c+1$ points in linear semi-uniform position. Then there exists a subset $\Gamma' \subset \Gamma$ of $(2c+1)$ points which spans $\P^c$ and is $3$-regular.
\end{proposition}

\begin{proof}
Since $\Gamma$ spans $\P^c$, we may assume that the $(c+1)$ coordinate points $p_0 , p_1 , \cdots , p_c$ of $\P^c$ are contained in $\Gamma$. We denote by $H_i$, $0 \leq i \leq c$, the hyperplane defined by $X_i$, or equivalently, the hyperplane spanned by $\{ p_0 , p_1 , \cdots , p_c \} - \{ p_i \}$.

If $\nu (c-1)=c$ and so $\Gamma$ is in linearly general position, then the result is well-known and can be easily verified.

Suppose that $\nu (i-1)=i$ and $\nu (i) \geq i+2$ for some $1 \leq i \leq c-1$. For $1 \leq j \leq i+1$, let
\begin{equation*}
\Delta_j = \langle \{p_0 , p_1 , \cdots , p_{i+1} \} - \{p_j \}\rangle \cong \P^i.
\end{equation*}
Since $\nu (i) \geq i+2$, there exists a point
\begin{equation*}
q_j = [a_{j,0} : a_{j,1} : \cdots : a_{j,c} ] \in \Gamma \cap \Delta_j - \{p_0 , p_1 , \cdots , p_c \}
\end{equation*}
for each $1 \leq j \leq i+1$. Moreover, it holds that $a_{j,k}=0$ if and only if $k \in \{j,i+2 , \cdots , c \}$ because $\nu (i-1)=i$. Also for each $i+1 \leq j \leq c$, there exists a point
\begin{equation*}
q_j = [a_{j,0} : a_{j,1} : \cdots : a_{j,c} ] \in \Gamma \cap \langle p_1 , \cdots , p_i , p_j \rangle  - \{p_0 , p_1 , \cdots , p_c \}
\end{equation*}
such that $a_{j,k}=0$ if and only if $k \in \{0,i+1,i+2,\cdots , c \} - \{j \}$. Now, let
\begin{equation*}
\Gamma' = \{p_0 , p_1 , \cdots , p_c \} \cup \{ q_1 , \cdots , q_c \}.
\end{equation*}
From the homogeneous coordinates of $q_j$'s, one can easily see that $\Gamma'$ consists of $(2c+1)$ distinct points and spans $\P^c$. From now on, let us prove that $\Gamma'$ is $3$-regular, or equivalently, that for each point $p \in \Gamma'$, there exists a quadratic hypersurface $Q_p \subset \P^c$ such that $Q_p \cap \Gamma' = \Gamma' - \{p \}$. Obviously, one can find such a quadratic hypersurface $Q_p$ if there exist two proper linear subspaces $\Lambda_1$ and $\Lambda_2$ of $\P^c$ such that $\Gamma'  - \{p \} \subset \Lambda_1 \cup \Lambda_2$ but $p \notin \Lambda_1 \cup \Lambda_2$.

For $i+2 \leq j \leq c$, observe that
\begin{equation*}
\Gamma' \cap H_j = \Gamma' - \{p_j , q_j \}.
\end{equation*}
Therefore we can find $Q_p$ for all $p \in \{ p_{i+2} , \cdots , p_c , q_{i+2} , \cdots , q_c \}$.

For $p=p_{i+1}$, observe that
\begin{equation*}
H_{i+1} \supset \Gamma' - \{p_{i+1} , q_1 , \cdots , q_i \}.
\end{equation*}
Also we get $p_{i+1} \notin \langle q_1 , \cdots , q_i \rangle$ since $\nu (i-1)=i$. Therefore we can take $\Lambda_1 = H_{i+1}$ and $\Lambda_2 = \langle q_1 , \cdots , q_i \rangle$.

For $p=p_j$ $(1 \leq j \leq i)$, note that $H_j$ contains $q_j$ and all $p_k$'s except $p_j$. Also letting \begin{equation*}
H_j ' := \langle ~ \{q_k ~|~ 1 \leq k \leq c,~ k \neq j \} ~ \rangle \quad (1 \leq j \leq i),
\end{equation*}
it can be easily checked that $p_j \notin H_j '$ by using the shape of the homogeneous coordinates of $q_k$'s $(1 \leq k \leq c)$. Therefore we can take $\Lambda_1 = H_j$ and $\Lambda_2 = H_j '$ for each $1 \leq j \leq i$.

For $p=q_j$ $(1 \leq j \leq i)$, note that $q_j \notin H_{j+1}$ and $q_{j+1} \in H_{j+1}$. Also it holds that
\begin{equation*}
q_j \notin H_j '' := \langle p_{j+1} , ~ \{ q_k ~|~ 1 \leq k \leq c , ~ k \neq j, j+1 \} ~ \rangle. 
\end{equation*}
Therefore we can take $\Lambda_1 = H_{j+1}$ and $\Lambda_2 = H_j ''$ for each $1 \leq j \leq i$.

For $p = q_{i+1}$, note that $q_{i+1} \notin H_1$ and $q_1 \in H_1$. Also it holds that
\begin{equation*}
q_{i+1} \notin H_{i+1} '' := \langle p_1 , ~ \{ q_k ~|~ 1 \leq k \leq c , ~ k \neq 1, i+1 \} ~ \rangle. 
\end{equation*}
Therefore we can take $\Lambda_1 = H_1$ and $\Lambda_2 = H_{i+1} ''$.

It remains to consider the case $p=p_0$. If $p_0 \notin \langle q_1 , \cdots , q_{i+1} \rangle$, then we can take $\Lambda_1 = H_0$ and $\Lambda_2 =\langle q_1 , \cdots , q_{i+1} \rangle$. Now, assume that $p_0 \in \langle q_1 , \cdots , q_{i+1} \rangle$. For each $1 \leq j \leq i+1$, we define
$$\Lambda_1 ^j := \begin{cases} \langle ~ \{ p_k ~|~ 1 \leq k \leq c, ~k \neq j+1 \}  , q_j \rangle & \mbox{if $j \neq i+1$, and}\\
                                \langle  \{ p_k ~|~ 1 \leq k \leq c, ~k \neq 1 \} , q_{i+1} \rangle & \mbox{if $j=i+1$} \end{cases}$$
and
$$\Lambda_2 ^j := \begin{cases} \langle H_j ' , p_{j+1} \rangle & \mbox{if $j \neq i+1$, and}\\
                                \langle H_{i+1} ' , p_1 \rangle & \mbox{if $j=i+1$.} \end{cases}$$
Then it is obvious that for all $1 \leq j \leq i+1$,
\begin{equation*}
p_0 \notin \Lambda_1 ^j
\end{equation*}
and
\begin{equation*}
\Gamma' - \{p_0 \} \subset \Lambda_1 ^j \cup \Lambda_2 ^j.
\end{equation*}
Thus it suffices to prove that $p_0 \notin \Lambda_2 ^j$ for some $1 \leq j \leq i+1$. Observe that $p_0 \in \Lambda_2 ^j$ implies $p_0 \in \Pi^j$ where
$$\Pi ^j := \begin{cases} \langle ~ \{ q_k ~|~ 1 \leq k \leq i+1, ~ k \neq j \} , p_{j+1} \rangle & \mbox{if $j \neq i+1$, and}\\
                                \langle ~ \{ q_k ~|~ 1 \leq k \leq i+1, ~ k \neq i+1 \}  , p_1 \rangle & \mbox{if $j=i+1$.} \end{cases}$$
Also $p_0 \in \Pi^j$ implies
\begin{equation*}
p_{j+1} \in \langle q_1 , \cdots , q_{i+1} \rangle \quad \mbox{if $1 \leq j \leq i$}
\end{equation*}
and
\begin{equation*}
p_1 \in \langle q_1 , \cdots , q_{i+1} \rangle \quad \mbox{if $j=i+1$}
\end{equation*}
since $p_0 \in \langle q_1 , \cdots , q_{i+1} \rangle$ and $\nu(i-1)=i$. Therefore if
\begin{equation*}
p_0 \in \bigcap_{j=1} ^{i+1} ~ \Lambda_2 ^j,
\end{equation*}
then we have $p_0 , p_1 , \cdots , p_{i+1} \in \langle q_1 , \cdots , q_{i+1} \rangle$, which is impossible. This completes the proof that $p_0 \notin \Lambda_2 ^j$ for at least one $j \in \{ 1, \cdots , i+1 \}$.
\end{proof}

\begin{remark}\label{rem:smallnumbercase}
Let $\Gamma \subset \P^c$ be a finite set of $d$ points in linear semi-uniform position. In the proof of Proposition \ref{prop:semiuniform}, it is shown that if $\nu (i) > i+1$ for some $1 \leq i \leq c-1$ then $d \geq 2c+1$. Therefore if $d \leq 2c$, then the linear semi-uniform position property implies that $\Gamma$ is in linearly general position and hence it is $3$-regular.
\end{remark}

\begin{corollary}\label{cor:degreeestimation}
Let $\Gamma \subset \P^c$ be a finite set of $d$ points in linear semi-uniform position. Then for each $m \geq 2$,
$$h^0 (\P^c , \mathcal{I}_{\Gamma} (m)) \begin{cases}
= {{c+m} \choose {m}}-d & \mbox{if $d \leq 2c+1$, and}\\
\leq {{c+m} \choose {m}}-2c-1 & \mbox{if $d \geq 2c+2$.} \end{cases}$$
In particular,
$$h^0 (\P^c , \mathcal{I}_{\Gamma} (2)) \begin{cases}
= {{c+1} \choose {2}}+c+1-d & \mbox{if $d \leq 2c+1$, and}\\
\leq {{c} \choose {2}} & \mbox{if $d \geq 2c+2$.} \end{cases}$$
\end{corollary}

\begin{proof}
If $d \leq 2c+1$, then $\Gamma$ is $3$-regular by Proposition \ref{prop:semiuniform} and Remark \ref{rem:smallnumbercase}. Therefore we get the desired equality. Now, suppose that $d \geq 2c+2$ and let $\Gamma'$ be as in Proposition \ref{prop:semiuniform}. Then we get
\begin{equation*}
h^0 (\P^c , \mathcal{I}_{\Gamma} (m)) \leq h^0 (\P^c , \mathcal{I}_{\Gamma'} (m)) = {{c+m} \choose {m}}-2c-1
\end{equation*}
where the latter equality follows immediately from the $3$-regularity of $\Gamma'$.
\end{proof}

We conclude this section by providing an example where Corollary \ref{cor:degreeestimation} is sharp.

\begin{example}
Let $\Gamma \subset \P^c$ be a general hyperplane section of a canonical curve $C \subset \P^{c+1}$ of genus $g=c+2$. Then $|\Gamma|=2c+2$ and $h^0 (\P^c , \mathcal{I}_{\Gamma}(2))={{c} \choose {2}}$. Thus, the inequality in Corollary
\ref{cor:degreeestimation} cannot be improved.
\end{example}
\vspace{0.5 cm}

\section{The deficiency module of curves}

\noindent Let $C \subset \P^{c+1}$ be a nondegenerate projective integral curve. The \textit{deficiency module} of $C$, denoted by $M(C)$, is defined by
\begin{equation*}
M(C) := \bigoplus_{m \geq 1} H^1 (\P^{c+1},\mathcal{I}_C (m)).
\end{equation*}
Also the Castelnuovo-Mumford regularity of $C$, denoted by ${\rm reg}(C)$, is the least integer $\ell$ such that $H^1 (\P^{c+1},\mathcal{I}_C (\ell-1))=H^1 (C,\mathcal{O}_C (\ell-2))=0$.

The aim of this section is to study the structure of $M(C)$ when the degree of $C$ is at most $2c$ and $C$ is not linearly normal and to apply it to Problem B in Section 1 for curve case.

We begin with recalling the following well-known fact.

\begin{proposition}\label{prop:lowdegcurve}
Let $C \subset \P^{c+1}$ be a nondegenerate integral projective curve of arithmetic genus $g$ and degree $d \leq 2c+1$. Then

\smallskip

\renewcommand{\descriptionlabel}[1]%
             {\hspace{\labelsep}\textrm{#1}}
\begin{description}
\setlength{\labelwidth}{13mm}
\setlength{\labelsep}{1.5mm}
\setlength{\itemindent}{0mm}

\item[(1)] $H^1 (C,\mathcal{O}_C (m))=0$ for all $m \geq 1$.
\item[(2)] If $C$ is linearly normal, then $M(C)=0$ and hence
\begin{equation}\label{eq:3.1}
a_m (C) = {{c+1+m} \choose {m}} - (md+1-g) \quad \mbox{for all $m \geq 2$.}
\end{equation}
\end{description}
\end{proposition}

\begin{proof}
(1) This follows from Clifford Theorem for projective integral curves (cf. \cite[Lemma 3.1]{KM}). \\
(2) By (1) and Riemann-Roch Theorem, we get $d=c+g+1$. Since $d \leq 2c+1$, it follows that $g \leq c$ and so $d \geq 2g+1$. Therefore $M(C)=0$ by \cite{C}(for smooth curves) and \cite{F}(for arbitrary integral curves). Now, (\ref{eq:3.1}) comes immediately by (1) and Riemann-Roch Theorem.
\end{proof}

By Proposition \ref{prop:lowdegcurve}.(2), $C$ is non-linearly normal if and only if $M(C)$ is nonzero. In such a case, ${\rm reg}(C) \geq 3$ and it is not hard to see that the sequence $h^1 (\P^{c+1},\mathcal{I}_C (m))$ decreases for all $m \geq 1$. The following theorem shows that one can say even more if $d \leq 2c$.

\begin{theorem}\label{thm:monotonousfunction}
Let $C \subset \P^{c+1}$ be a nondegenerate integral projective curve of arithmetic genus $g$ and degree $d \leq 2c$. Suppose that $C$ is not linearly normal. Then

\smallskip

\renewcommand{\descriptionlabel}[1]%
             {\hspace{\labelsep}\textrm{#1}}
\begin{description}
\setlength{\labelwidth}{13mm}
\setlength{\labelsep}{1.5mm}
\setlength{\itemindent}{0mm}

\item[(1)] For all $2 \leq m \leq {\rm reg}(C)-1$,
\begin{equation*}
h^1 (\P^{c+1},\mathcal{I}_C (m-1)) > h^1 (\P^{c+1},\mathcal{I}_C (m)).
\end{equation*}
\item[(2)] For all $2 \leq m \leq {\rm reg}(C)-1$,
\begin{equation*}
{\rm reg}(C) \leq m + h^1 (\P^{c+1},\mathcal{I}_C (m-1)).
\end{equation*}
\item[(3)] ${\rm reg}(C) \leq d-c+1-g$. Moreover, ${\rm reg}(C) = d-c+1-g$ if and only if
\begin{equation*}
h^1 (\P^{c+1},\mathcal{I}_C (m)) = d-c-g-m \quad \mbox{for all $1 \leq m \leq d-c-g$.}
\end{equation*}
\end{description}
\end{theorem}
\smallskip

\begin{proof}
Let $\Gamma \subset \P^c$ be a general hyperplane section of $C$ defined by a linear form $H$ on $\P^{c+1}$. By Remark \ref{rem:smallnumbercase}, $\Gamma$ is in linearly general position and so ${\rm reg}(\Gamma)\leq 3$. Moreover, its homogeneous ideal $I_{\Gamma}$ is generated by quadrics since $d \leq 2c$ (cf. \cite[Theorem 1]{GL}). Let $R$ be the homogeneous coordinate ring of $\P^{c+1}$ and consider the graded $R$-module
\begin{equation*}
E := I_{\Gamma} / \langle I_C , H \rangle.
\end{equation*}
If we apply sheaf cohomology to the short exact sequence
\begin{equation*}
0 \rightarrow \mathcal{I}_C (-1) \rightarrow \mathcal{I}_C \rightarrow \mathcal{I}_{\Gamma} \rightarrow 0
\end{equation*}
we obtain a grade-preserving exact sequence
\begin{equation*}
0 \rightarrow E \rightarrow M(C)(-1) \rightarrow  M(C) \rightarrow 0
\end{equation*}
since $\Gamma$ is $3$-regular. In particular, it holds that for each $m \geq 2$,
\begin{equation*}
h^1 (\P^{c+1},\mathcal{I}_C (m-1)) \geq h^1 (\P^{c+1},\mathcal{I}_C (m))
\end{equation*}
and the inequality turns into the equality if and only if $E_m =0$. Note that $E$ is generated by homogeneous elements of degree $2$ since $I_{\Gamma}$ is generated by quadrics. Thus if $E_m = 0$ for some $m \geq 2$, then $E_{m'} =0$ for all $m' \geq m$ and hence
\begin{equation*}
h^1 (\P^{c+1},\mathcal{I}_C (m-1)) = h^1 (\P^{c+1},\mathcal{I}_C (m')) \quad \mbox{for all $m' \geq m$}
\end{equation*}
Obviously this occurs if and only if $h^1 (\P^{c+1},\mathcal{I}_C (m-1))=0$. Consequently, it is shown that if $h^1 (\P^{c+1},\mathcal{I}_C (m-1))\neq 0$, or equivalently (by Proposition \ref{prop:lowdegcurve}.(1)), if $2 \leq m \leq {\rm reg}(C)-1$, then
\begin{equation*}
h^1 (\P^{c+1},\mathcal{I}_C (m-1)) > h^1 (\P^{c+1},\mathcal{I}_C (m)).
\end{equation*}
Now, let us consider the sequence of positive integers
\begin{equation*}
b_m (C) := h^1 (\P^{c+1},\mathcal{I}_C (m-1)) - h^1 (\P^{c+1},\mathcal{I}_C (m))
\end{equation*}
for $2 \leq m \leq {\rm reg}(C)-1$. Then it holds that
\begin{equation*}
{\rm reg}(C) \leq m + \sum_{j=m} ^{{\rm reg}(C)-1} b_m (C) = m + h^1 (\P^{c+1},\mathcal{I}_C (m-1))
\end{equation*}
for all $2 \leq m \leq {\rm reg}(C)-1$. In particular,
\begin{equation*}
{\rm reg}(C) \leq 2+ h^1 (\P^{c+1},\mathcal{I}_C (1))=d-c+1-g
\end{equation*}
and the inequality turns into the equality if and only if $b_m (C)=1$ for all $2 \leq m \leq {\rm reg}(C)-1$. This completes the proof of (3).
\end{proof}

\begin{example}
For the rational normal surface scroll $S:= S(1,c-1) \subset \P^{c+1}$, let $H$ and $F$ be respectively the hyperplane section and a ruling of $S$ and let $C$ be an irreducible divisor of $S$ linearly equivalent to $H+(c+1)F$. Then
\begin{equation*}
C \subset \P^{c+1}
\end{equation*}
is a smooth rational curve of degree $2c+1$ such that
$$h^1 (\P^{c+1},\mathcal{I}_C (m)) = \begin{cases} c & \mbox{if $m=1$ or $2$,}\\
                                                        c-m+1 & \mbox{if $3 \leq m \leq c$, and}\\
                                                        0 & \mbox{if $m >c$.}\end{cases}$$
In particular, $h^1 (\P^{c+1},\mathcal{I}_C (1)) = h^1 (\P^{c+1},\mathcal{I}_C (2))$. Thus, the hypothesis $d \leq 2c$ in Theorem \ref{thm:monotonousfunction} cannot be weakened.
\end{example}

\begin{remark}
(1) Theorem \ref{thm:monotonousfunction}.(1) and the second part of Theorem \ref{thm:monotonousfunction}.(3) are known to be true when $g=0$. See \cite[Theorem 3.3 and Proposition 3.5]{BS2}

\noindent (2) A part of A. Noma's result in \cite{No} says that the inequality
\begin{equation*}
{\rm reg}(C) \leq d-c+1-g
\end{equation*}
holds under the condition $g \leq c-1$ which is satisfied if $d \leq 2c$. So, the first part of Theorem \ref{thm:monotonousfunction}.(3) is a reproof of A. Noma's result for the case $d \leq 2c$. Note that our proof is obtained by using the hyperplane section instead of the vector bundle technique. The second part of Theorem \ref{thm:monotonousfunction}.(3) provides a cohomological characterization of the extremal cases with respect to the above inequality.
\end{remark}

Theorem \ref{thm:monotonousfunction} enables us to obtain the following satisfactory answer for Problem B in Section 1 for curve case.

\begin{theorem}\label{thm:curvehomologicaltype}
Let $c$ and $k$ be integers such that $c \geq 2$ and $1 \leq k \leq c$. Then

\smallskip

\renewcommand{\descriptionlabel}[1]%
             {\hspace{\labelsep}\textrm{#1}}
\begin{description}
\setlength{\labelwidth}{13mm}
\setlength{\labelsep}{1.5mm}
\setlength{\itemindent}{0mm}

\item[(1)] $\delta_{1,c,2} (k) = {{c+1} \choose {2}} +1-k$.
\item[(2)] Let $C \subset \P^{c+1}$ be a nondegenerate projective integral curve of arithmetic genus $g$ and degree $d$. Then
    \begin{equation*}
    a_2 (C) = \delta_{1,c,2} (k) ~ \Longleftrightarrow ~ h^1 (\P^{c+1},\mathcal{I}_C (2))=2(d-c)-1-g-k.
    \end{equation*}
    Furthermore, if $a_2 (C) = \delta_{1,c,2} (k)$ holds, then
    \begin{enumerate}
    \item[(a)] $C$ is linearly normal if and only if $g=k-1$ and $d=c+k$; and
    \item[(b)] if $C$ is not linearly normal then
    \begin{equation*}
     0 \leq g \leq k-3 \quad \mbox{and} \quad c+\frac{g+k+1}{2} \leq d < c+k.
    \end{equation*}
    \end{enumerate}
\end{description}
\end{theorem}
\smallskip

\begin{proof}
(1) We get
\begin{equation*}
\delta_{1,c,2} (1) \leq {{c+1} \choose {2}}
\end{equation*}
from (\ref{eq:1.1}). Also for any linearly normal curve $C \subset \P^{c+1}$ of arithmetic genus $k-1$ and degree $c+k$, it holds that
\begin{equation*}
a_2 (C) = {{c+1} \choose {2}} +1-k
\end{equation*}
(eg. Proposition \ref{prop:lowdegcurve}.(2)). That is,
\begin{equation*}
{{c+1} \choose {2}} +1-k \in A_{1,c,2} \quad \mbox{for all $1 \leq k \leq c$.}
\end{equation*}
Obviously, this completes the proof.

\noindent (2) Let $\Gamma \subset \P^c$ be a general hyperplane section of $C$. If $d \geq c+k+1$, then
\begin{equation*}
a_2 (C) \leq h^0 (\P^c , \mathcal{I}_C (2)) \leq {{c+1} \choose {2}} - (c+k+1) < \delta_{1,c,2} (k).
\end{equation*}
Thus the condition $a_2 (C) = \delta_{1,c,2}(k)$ implies that $d \leq c+k \leq 2c$. Then we get
\begin{equation}\label{eq:3.2}
h^1 (\P^{c+1},\mathcal{I}_C (1))=d-g-c-1
\end{equation}
and
\begin{equation}\label{eq:3.3}
a_2 (C) = {{c+1} \choose {2}}-\{2d-2c-2-g-h^1 (\P^{c+1},\mathcal{I}_C (2)) \}
\end{equation}
(cf. Proposition \ref{prop:lowdegcurve}.(1) and Riemann-Roch Theorem). Now, the first assertion comes by comparing $\delta_{1,c,2}(k)$ with $a_2 (C)$ in (\ref{eq:3.3}).

(2.a): If $C$ is linearly normal, then we get $g=k-1$ by comparing $\delta_{1,c,2} (k)$ with $a_2 (C)$ in (\ref{eq:3.1}). Conversely, if $g=k-1$ and $d=c+k$ then $C$ is linearly normal by (\ref{eq:3.2}).
\smallskip

(2.b): Suppose that $C$ is not linearly normal. The inequality $c+ \frac{g+k+1}{2} \leq d$ is obtained from
\begin{equation*}
h^1 (\P^{c+1},\mathcal{I}_C (2))=2(d-c)-1-g-k \geq 0.
\end{equation*}
Also it holds from (\ref{eq:3.2}) that $g < d-c-1$. By Theorem \ref{thm:monotonousfunction}, we have the inequality
\begin{equation*}
h^1 (\P^{c+1},\mathcal{I}_C (1)) > h^1 (\P^{c+1},\mathcal{I}_C (2)),
\end{equation*}
which implies $d < c+k$ (cf. (\ref{eq:3.2}) and (\ref{eq:3.3})). These complete the proof.
\end{proof}

Theorem \ref{thm:curvehomologicaltype} says that if there exists a non-linearly normal projective integral curve $C \subset \P^{c+1}$ of degree $d$ and arithmetic genus $g$ such that $a_2 (C)= \delta_{1,c,2} (k)$ for some $1 \leq k \leq c$, then $k \geq 3$ and the pair $(g,d)$ should be contained in the shadowed region of Figure 1. In such a case, we will say that $(g,d)$ is \textit{realizable}. Along this line, the remaining part of this section is devoted to consider two problems. Firstly, one can naturally ask if a given pair $(g,d)$ of integers in the shadowed region of Figure 1 is realizable or not. This problem seems to be very hard in general (cf. Remark \ref{rm:difficulty}). In this direction, we will show that all pairs are realizable if $3 \leq k \leq 7$ or if $(g,d)$ lies on the lines $g=0$ or $d=c+k-1$. Secondly, we will provide a geometric description of all curves $C \subset \P^{c+1}$ with $a_2 (C) = \delta_{1,c,2}(k)$ when $1 \leq k \leq 4$ (see Theorem \ref{thm:quadric-curve}).

\begin{figure}[t]
\begin{center}
\begin{pspicture}(-1,0)(11,8)
\footnotesize \psset{xunit=.5cm,yunit=.5cm}

\psaxes[labels=none,ticks=none]{->}(0,0)(-1,-1)(20,14)

\uput[180](-.5,4){$c+1+\frac{k}{2}$}
\uput[180](-.5,10.6){$c+k$}
\uput[180](-.5,14){$d$}
\uput[180](12,5){$d=\frac{1}{2} g + (c+1+\frac{k}{2})$}

\rput[t](.5,-.5){0} \rput[t](12,-.5){$k-2$} \rput(21,0){$g$}

\psline[linestyle=dashed](11.9,0)(11.9,6) \psline(11.9,6)(11.9,14)
\pspolygon[linestyle=none,fillstyle=hlines,hatchsep=6pt](0,10.6)(0,4)(11.9,10.6)

\psline(0,4)(18,14)
\psline(0,10.6)(18,10.6)
\psline[arrows=->](7,6)(6,7)

\end{pspicture}
\end{center}
\caption{Possible pairs of $(g,d)$ for a given $(c,k)$}
\end{figure}
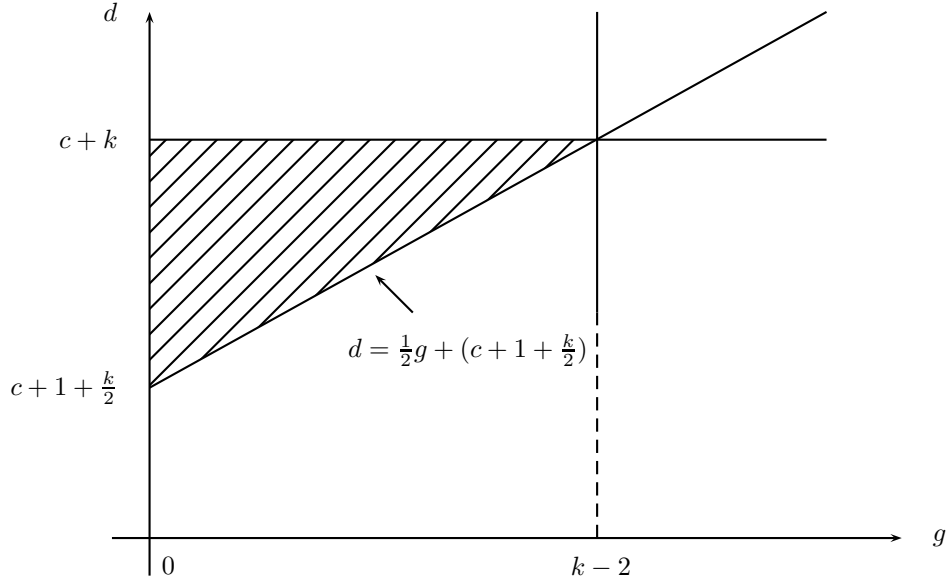

The following two examples show that any $(g,d)$ on the lines $g=0$ and $d=c+k-1$ are realizable.

\begin{example}\label{ex:genuszero}
Let $k$ and $d$ be integers satisfying $3 \leq k \leq c$ and $c+\frac{k+1}{2} \leq d < c+k$. For the smooth rational normal surface scroll
\begin{equation*}
S:=S(c+k-d,d-k) \subset \P^{c+1},
\end{equation*}
let $C$ be an irreducible divisor on $S$ linearly equivalent to $H+(d-c)F$ where $H$ and $F$ be respectively a hyperplane section and a ruling of $S$. Then $C \subset \P^{c+1}$ is a nondegenerate smooth rational curve of degree $d$ such that $a_2 (C) = \delta_{1,c,2}(k)$.
\end{example}

\begin{example}\label{extremalsecantline}
Let $k$ and $g$ be integers satisfying $1 \leq k \leq c$ and $0 \leq g \leq k-3$. Note that $d:=c+k-1$ is at least $2g+5$ and hence there exists a linearly normal smooth curve
\begin{equation*}
\widetilde{C} \subset \P^{c+k-1-g}
\end{equation*}
of genus $g$ and degree $d$. Choose general $(k-g)$-points $p_1 , \cdots , p_{k-g}$ on $\widetilde{C}$. Then $\Delta := \langle p_1 , \cdots , p_{k-g} \rangle \subset \P^{c+k-1-g}$ is a $(k-1-g)$-dimensional subspace. A general $(k-3-g)$-dimensional subspace $\Lambda \subset \Delta$ defines an isomorphic projection of $\widetilde{C}$ to a curve $C:= \pi_{\Lambda} (\widetilde{C}) \subset \P^{c+1}$ of genus $g$ and degree $d$.
Furthermore, $\pi_{\Lambda} (\Delta \setminus \Lambda)$ is a $(k-g)$-secant line to $C$ and hence the regularity of $C$ is at least $k-g=d-c+1-g$. On the other hand, Theorem \ref{thm:monotonousfunction}.(3) shows that ${\rm reg}(C) \leq d-c+1-g$ since $d \leq 2c$. Thus, we get ${\rm reg}(C)= d-c+1-g$. Moreover, it holds that
\begin{equation*}
h^1 (\P^{c+1},\mathcal{I}_C (m))=d-c-g-m \quad \mbox{for all $1 \leq m \leq d-c-g$}
\end{equation*}
again by Theorem \ref{thm:monotonousfunction}.(3). In particular, $a_2 (C) = \delta_{1,c,2} (k)$ by Theorem \ref{thm:curvehomologicaltype}.(2).
\end{example}

Concerned with the first problem for $3 \leq k \leq 7$, let us consider Table 1 which is obtained by Theorem \ref{thm:curvehomologicaltype}.(2).\\

\begin{table}[hbt]
\begin{center}
\begin{tabular}{|c|c|c|c|c| }\hline
$\quad k \quad$ & $\quad g \quad$ & $\quad d \quad$ &  $h^1 (\P^{c+1},\mathcal{I}_C (1))$ & $h^1 (\P^{c+1},\mathcal{I}_C (2))$   \\\hline
$3$             & $0$             &       $c+2$     &                  $1$                &   $0$        \\\hline
$4$             & $0$             &       $c+3$     &                  $2$                &   $1$        \\\cline{2-5}
                & $1$             &       $c+3$     &                  $1$                &   $0$        \\\hline
                & $0$             &       $c+3$     &                  $2$                &   $0$        \\\cline{3-5}
$5$             &                 &       $c+4$     &                  $3$                &   $2$        \\\cline{2-5}
                & $1$             &       $c+4$     &                  $2$                &   $1$        \\\cline{2-5}
                & $2$             &       $c+4$     &                  $1$                &   $0$        \\\hline
                & $0$             &       $c+4$     &                  $3$                &   $1$        \\\cline{3-5}
                &                 &       $c+5$     &                  $4$                &   $3$        \\\cline{2-5}
$6$             & $1$             &       $c+4$     &                  $2$                &   $0$       \\\cline{3-5}
                &                 &       $c+5$     &                  $3$                &   $2$       \\\cline{2-5}
                & $2$             &       $c+5$     &                  $2$                &   $1$       \\\cline{2-5}
                & $3$             &       $c+5$     &                  $1$                &   $0$       \\\hline
                & $0$             &       $c+4$     &                  $3$                &   $0$        \\\cline{3-5}
                &                 &       $c+5$     &                  $4$                &   $2$        \\\cline{3-5}
                &                 &       $c+6$     &                  $5$                &   $4$        \\\cline{2-5}
$7$             & $1$             &       $c+5$     &                  $3$                &   $1$       \\\cline{3-5}
                &                 &       $c+6$     &                  $4$                &   $3$       \\\cline{2-5}
                & $2$             &       $c+5$     &                  $2$                &   $0$       \\\cline{3-5}
                &                 &       $c+6$     &                  $3$                &   $2$       \\\cline{2-5}
                & $3$             &       $c+6$     &                  $2$                &   $1$       \\\cline{2-5}
                & $4$             &       $c+6$     &                  $1$                &   $0$       \\\hline
\end{tabular}
\end{center}
\caption{Non-linearly normal curve with $a_2 (C) = \delta_{1,c,2}(k)$}
\end{table}

\noindent By Example \ref{ex:genuszero} and Example \ref{extremalsecantline}, all pairs $(g,d)$ in Table 1 are realizable except the following three cases:
$$\begin{cases} (i) & k=6 \quad \mbox{and} \quad (g,d)=(1,c+4);\\
                (ii) & k=7 \quad \mbox{and} \quad (g,d)=(1,c+5);\\
                (iii) & k=7 \quad \mbox{and} \quad (g,d)=(2,c+5). \end{cases}$$
For these cases, we provide the following examples.

\begin{example}\label{ex:genusone}
Suppose that $c \geq 5$. Let $E$ be an elliptic curve and let $\mathcal{L}_1$ and $\mathcal{L}_2$ be line bundles on $E$ of degree $3$ and $c-1$, respectively. Consider the elliptic normal surface scroll
\begin{equation*}
S:= \P_E (\mathcal{L}_1 \oplus \mathcal{L}_2 ) \subset \P^{c+1}.
\end{equation*}
Note that $S$ is projectively normal. Let $H$ and $F$ be respectively a hyperplane section and a ruling of $S$ and consider an irreducible divisor $C$ on $S$ linearly equivalent to $H+ \mathcal{L} F$ where $\mathcal{L}$ is a line bundle on $E$ of degree $\ell \geq 2$. Then $C \subset \P^{c+1}$ is a nondegenerate elliptic curve of degree $d:=c+2+\ell$ and
\begin{equation*}
h^1 (\P^{c+1},\mathcal{I}_C (2))= h^1 (C,\mathcal{L}_1 \otimes \mathcal{L}^{-1} ) +  h^1 (C,\mathcal{L}_2 \otimes \mathcal{L}^{-1} ).
\end{equation*}
If $\ell=2$, then $h^1 (\P^{c+1},\mathcal{I}_C (2))=0$ and so $a_2 (C) = \delta_{1,c,2}(6)$. Also if $\mathcal{L}=\mathcal{L}_1$, then $h^1 (\P^{c+1},\mathcal{I}_C (2))=1$ and so $a_2 (C) = \delta_{1,c,2}(7)$.
Therefore, the cases $(i)$ and $(ii)$ are realizable.
\end{example}

\begin{example}\label{ex:genusone}
Suppose that $c \geq 5$. Let $\mathcal{C}$ be a smooth curve of genus $2$ and let $\mathcal{L}_1$ and $\mathcal{L}_2$ be line bundles on $\mathcal{C}$ of degree $5$ and $c$, respectively. Consider the linearly normal surface scroll
\begin{equation*}
S:= \P_{\mathcal{C}} (\mathcal{L}_1 \oplus \mathcal{L}_2 ) \subset \P^{c+2}.
\end{equation*}
and its general hyperplane section $C$. Thus $C \subset \P^{c+1}$ is a smooth curve of genus $2$ and degree $c+5$. Since $S$ is projectively normal and $3$-regular, it holds that $h^1 (\P^{c+1},\mathcal{I}_C (2))=0$ and so $a_2 (C)= \delta_{1,c,2}(7)$. Therefore, case $(iii)$ is realizable.
\end{example}

\begin{remark}\label{rm:difficulty}
For $k=8$, the pair $(g,d)=(1,c+6)$ is contained in the shadowed region of Figure 1. This pair is realizable if and only if there exists a nondegenerate projective integral curve $C \subset \P^{c+1}$ of arithmetic genus $1$ and degree $c+6$ such that $h^1 (\P^{c+1},\mathcal{I}_C (1))=4$ and $h^1 (\P^{c+1},\mathcal{I}_C (2))=2$. We don't know yet if such a curve exists or not.
\end{remark}

Now, we turn to the second problem outlined above.

\begin{theorem}\label{thm:quadric-curve}
Let $C \subset \P^{c+1}$ be a nondegenerate projective integral curve of degree $d$. Then

\smallskip

\renewcommand{\descriptionlabel}[1]%
             {\hspace{\labelsep}\textrm{#1}}
\begin{description}
\setlength{\labelwidth}{13mm}
\setlength{\labelsep}{1.5mm}
\setlength{\itemindent}{0mm}

\item[(1)] $a_2 (C) = {{c+1} \choose {2}}$ if and only if $C$ is a rational normal curve.
\item[(2)] $a_2 (C) = {{c+1} \choose {2}}-1$ if and only if $C$ is a linearly normal curve of arithmetic genus one.
\item[(3)] Suppose that $c \geq 3$. Then $a_2 (C) = {{c+1} \choose {2}}-2$ if and only if $C$ is either
   \begin{enumerate}
    \item[(i)] the image of an isomorphic projection of a rational normal curve from a point
    \end{enumerate}
    or else
    \begin{enumerate}
    \item[(ii)] a linearly normal curve of arithmetic genus two.
    \end{enumerate}
\item[(4)] Suppose that $c \geq 4$. Then $a_2 (C) = {{c+1} \choose {2}}-3$ if and only if $C$ is either
    \begin{enumerate}
    \item[(i)] a smooth rational curve of $d=c+3$ having a $4$-secant line
    \end{enumerate}
    or
    \begin{enumerate}
    \item[(ii)] the image of an isomorphic projection of a linearly normal curve of arithmetic genus one from a point
    \end{enumerate}
    or else
    \begin{enumerate}
    \item[(iii)] a linearly normal curve of arithmetic genus three.
    \end{enumerate}
\end{description}
\end{theorem}

\begin{proof}
All the statements are obtained by combining Theorem \ref{thm:curvehomologicaltype} and the geometric reinterpretation of cohomological types listed in Table 1. For the case where $(k,g,d)=(4,0,c+3)$, we refer the reader to \cite[Corollary 4.3]{BS1}.
\end{proof}
 
\begin{remark}
(1) Let $C \subset \P^{c+1}$ be a non-linearly normal curve. Since $C$ is obtained as the image of an isomorphic linear projection of a linearly normal curve, it is natural to understand $C$ by the location of the projection center. For example, $C$ has the triple $(g,d,h^1 (\P^{c+1},\mathcal{I}_C (2)))=(0,c+3,2)$ if and only if it is the image of an isomorphic projection of a rational normal curve $\widetilde{C} \subset \P^{c+3}$ from a line $L$ contained in a $4$-secant $3$-space to $\widetilde{C}$ (cf. \cite[Corollary 4.3.(c)]{BS1}).

\noindent (2) In order to extend Theorem \ref{thm:quadric-curve} to the next case, we need to consider four cases in Table 1. For the case where $(g,d,h^1 (\P^{c+1},\mathcal{I}_C (2)))=(0,c+3,2)$ or $(2,c+4,0)$, this is easy. But we don't know yet a precise answer for the other two cases.
\end{remark}
\vspace{0.5 cm}

\section{Proof of Theorem \ref{thm:main}}

\noindent The aim of this section is to give a proof of Theorem \ref{thm:main}.

\begin{definition and remark}\label{lowdegreevarieties}
Let $X \subset \P^{n+c}$ be a nondegenerate projective variety of degree $d$ and codimension $c \geq 2$. It is well-known that $d \geq c+1$.
\smallskip

\noindent (A) $X$ is called a \textit{variety of minimal degree} if $d=c+1$. A variety of minimal degree is either (a cone over) the Veronese surface in $\P^5$ or a rational normal scroll. If $X$ is a variety of minimal degree, then it is arithmetically Cohen-Macaulay and
\begin{equation*}
a_m (X) = F(n,c,m) \quad \mbox{for all $m \geq 1$.}
\end{equation*}
Note that
\begin{equation}\label{eq:4.1}
F(n,c,2) = {{c+1} \choose {2}}
\end{equation}
and
\begin{equation}\label{eq:4.2}
F(n,c,m)=F(n,c,m-1)+F(n-1,c,m) \quad \mbox{for all $m \geq 2$.}
\end{equation}

\noindent (B) $X$ is called a \textit{variety of almost minimal degree} if $d=c+2$. Also $X$ is called a \textit{del Pezzo variety} if it is
a variety of almost minimal degree such that ${\rm depth}(X)=n+1$. We refer to \cite{BP} and \cite{F2} for the classification of varieties of almost minimal degree. If $X$ is a variety of almost minimal degree, then
\begin{equation*}
a_m (X) = G_t (n,c,m) \quad \mbox{for all $m \geq 2$}
\end{equation*}
where $t$ denotes the arithmetic depth of $X$ (cf. \cite[Theorem A and B]{HSV} or \cite[Theorem 2.2]{N}). Note that
\begin{equation}\label{eq:4.3}
G_t (n,c,2) = {{c+1} \choose {2}}+t-n-2
\end{equation}
and
\begin{equation}\label{eq:4.4}
G_t (n,c,m)=G_t (n,c,m-1)+G_t (n-1,c,m) \quad \mbox{for all $m \geq 2$.}
\end{equation}

\noindent (C) For each $m \geq 2$, it holds that
\begin{equation}\label{eq:4.5}
F(n,c,m) > G_{n+1} (n,c,m) > G_n (n,c,m) > \cdots > G_1 (n,c,m).
\end{equation} 

\noindent (D) For each $1 \leq k \leq c+1$, we define
\begin{equation*}
H_k (n,c,m) = {{m+n+c} \choose {m}} \quad \quad \quad \quad \quad \quad \quad \quad \quad \quad \quad \quad \quad \quad \quad \quad \quad \quad
\end{equation*}
\begin{equation*}
\quad \quad \quad \quad \quad  - \{ (c+k) {{m+n-1} \choose {n}} + (2-k){{m+n-2} \choose {n-1}} +{{m+n -2} \choose {n-2}} \}.
\end{equation*}
Observe that
\begin{equation}\label{eq:4.6}
H_k (n,c,2) = {{c+1} \choose {2}}+1-k
\end{equation}
and
\begin{equation}\label{eq:4.7}
H_k (n,c,m)=H_k (n,c,m-1)+ H_k (n-1,c,m) \quad \mbox{for all $m \geq 2$.}
\end{equation}
It is probably well-known that if $X$ is arithmetically Cohen-Macaulay and $d=c+k$, then
\begin{equation}\label{eq:4.8}
a_m (X) = H_k (n,c,m) \quad \mbox{for all $m \geq 1$}.
\end{equation}
In particular, 
\begin{equation}\label{eq:4.9}
F (n,c,m) = H_1 (n,c,m) \quad \mbox{and} \quad G_{n+1} (n,c,m) = H_2 (n,c,m).
\end{equation}
It is easy to check that
\begin{equation}\label{eq:4.10}
H_1 (n,c,m)> H_2 (n,c,m)> \cdots > H_c (n,c,m) \quad \mbox{for all $m \geq 2$.}
\end{equation}

\noindent (E) Observe that
\begin{equation*}
G_n (n,c,m) = H_3 (n,c,m) + {{m+n-3} \choose {n}}.
\end{equation*}
Therefore it holds that
\begin{equation}\label{eq:4.11}
G_n (n,c,2) = H_3 (n,c,2)
\end{equation}
and
\begin{equation}\label{eq:4.12}
G_n (n,c,m) > H_3 (n,c,m) \quad \mbox{for all $m \geq 3$.}
\end{equation}

\noindent (F) From (\ref{eq:4.5}), (\ref{eq:4.10}), (\ref{eq:4.11}) and (\ref{eq:4.12}), we get the following for all $m \geq 2$:
\begin{equation}\label{eq:4.13}
\begin{CD}
H_1 (n,c,m)      & =    & F(n,c,m)              &  \\
\vee             &      &  \vee                 & \\
H_2 (n,c,m)      & =    & \quad G_{n+1} (n,c,m) & \\
\vee             &      &  \vee                 & \\
H_3 (n,c,m)      & \leq &  G_n (n,c,m)          &  \mbox{with equality only when $m=2$}\\
\vee             &      &   \vee                & \\
\vdots           &      &  \vdots               & \\
\vee             &      &  \vee                 & \\
H_{c+1} (n,c,m)  &      &  G_1 (n,c,m)          &
\end{CD}
\end{equation}
\end{definition and remark}
\vspace{0.3 cm}

Keeping Definition and Remark \ref{lowdegreevarieties} in mind, we begin with the following useful result.

\begin{proposition}\label{prop:degreebound}
Let $c$ and $k$ be integers such that $1 \leq k \leq c+1$. If $X \subset \P^{n+c}$ is a nondegenerate projective variety of codimension $c$ and degree $d \geq c+k$, then
\begin{equation*}
a_m (X) \leq H_k (n,c,m) \quad \mbox{for all $m \geq 2$.}
\end{equation*}
\end{proposition}

\begin{proof}
We will prove our theorem by the induction on $m$ and $n$. Let $Y \subset \P^{n+c-1}$ be a general hyperplane section of $X$. Throughout the proof, we will use the basic inequality
\begin{equation*}
a_m (X) \leq a_{m-1} (X) + a_m (Y),
\end{equation*}
which follows from the exact sequence $0 \rightarrow \mathcal{I}_X (-1) \rightarrow \mathcal{I}_X \rightarrow \mathcal{I}_Y \rightarrow 0$ where $\mathcal{I}_Y$ is the sheaf of ideals of $Y$ in $\P^{n+c-1}$.

For $m=2$, let $\Lambda$ be a general $c$-dimensional linear subspace of $\P^{n+c}$ and consider the finite set
\begin{equation*}
\Gamma := X \cap \Lambda \subset \P^c
\end{equation*}
of $d$ points in linear semi-uniform position. By Corollary \ref{cor:degreeestimation}, we have
$$a_2 (X) \leq  a_2 (\Gamma) \begin{cases}
= {{c+1} \choose {2}}+c+1-d & \mbox{if $d \leq 2c+1$, and}\\
\leq {{c} \choose {2}} & \mbox{if $d \geq 2c+2$.} \end{cases}$$
This shows that
\begin{equation*}
a_2 (X) \leq H_k (n,c,2)= {{c+1} \choose {2}}+c+1-(c+k)
\end{equation*}
since $k \leq c+1$ and $d \geq c+k$ (cf. (\ref{eq:4.6})).

Suppose that $n=1$ and $m \geq 3$. Then
\begin{equation*}
\begin{CD}
a_m (X) & \quad \leq \quad & a_{m-1} (X) + a_m (Y) \\
        & \quad \leq \quad & H_k (1,c,m-1) +  {{m+c} \choose {m}} - (d+k) = H_k (1,c,m)
\end{CD}
\end{equation*}
by the induction hypothesis and Corollary \ref{cor:degreeestimation}.

Finally, suppose that $n \geq 2$ and $m \geq 3$. Then
\begin{equation*}
\begin{CD}
a_m (X) & \quad \leq \quad & a_{m-1} (X) + a_m (Y) \\
        & \quad \leq \quad & H_k (n,c,m-1) +  H_k (n-1,c,m) = H_k (n,c,m)
\end{CD}
\end{equation*}
by the induction hypothesis and (\ref{eq:4.7}).
\end{proof}

\begin{corollary}\label{cor:quadraticinequality}
Let $X \subset \P^{n+c}$ be a nondegenerate projective variety of codimension $c$ and degree $d$. If
\begin{equation*}
a_2 (X) = {{c+1} \choose {2}}+1-k \quad \mbox{for some $1 \leq k \leq c$,}
\end{equation*}
then $d \leq c+k$.
\end{corollary}

\begin{proof}
This follows immediately by Proposition \ref{prop:degreebound} (cf. (\ref{eq:4.6})).
\end{proof}

\begin{corollary}\label{cor:boundarycases}
Let $c$ and $k$ be integers such that $1 \leq k \leq c$ and let $X \subset \P^r$ be a nondegenerate projective variety of dimension $n$, codimension $c>0$ and degree $d=c+k$. Then
\begin{equation*}
a_2 (X) = {{c+1} \choose {2}}+1-k \quad \Longleftrightarrow \quad {\rm depth}(X)=n+1.
\end{equation*}
\end{corollary}

\begin{proof}
Let $C := X \cap \P^{c+1} \subset \P^{c+1}$ where $\P^{c+1}$ is a general linear subspace of $\P^r$ and let $\Gamma \subset \P^c$ be a general hyperplane section of $C$. Thus we have
\begin{equation}\label{eq:4.14}
a_2 (X) \leq a_2 (C) \leq h^0 (\P^c , \mathcal{I}_{\Gamma} (2)) = {{c+1} \choose {2}}+1-k
\end{equation}
(cf. Corollary \ref{cor:degreeestimation}). Note that ${\rm depth}(X)=n+1$ if and only if ${\rm depth}(C)=2$.
\smallskip

($\Longrightarrow$): By (\ref{eq:4.14}), we know that $a_2 (C)={{c+1} \choose {2}}+1-k$. If $C$ is not linearly normal, then $d < c+k$ by Theorem \ref{thm:curvehomologicaltype}.(2), a contradiction. Therefore $C$ is linearly normal and so ${\rm depth}(C)=2$ by Proposition \ref{prop:lowdegcurve}.(2).
\smallskip

($\Longleftarrow$): If ${\rm depth}(X)=n+1$, then the inequalities in (\ref{eq:4.14}) turn into equalities. This completes the proof.
\end{proof}

From Proposition \ref{prop:degreebound}, we obtain the following main result of this section.

\begin{theorem}\label{thm:firstandsecond}
Let $X \subset \P^{n+c}$ be a nondegenerate projective variety of degree $d$ and codimension $c \geq 2$, and let $m \geq 2$ be an integer. Then

\smallskip

\renewcommand{\descriptionlabel}[1]%
             {\hspace{\labelsep}\textrm{#1}}
\begin{description}
\setlength{\labelwidth}{13mm}
\setlength{\labelsep}{1.5mm}
\setlength{\itemindent}{0mm}

\item[(1)] $a_m (X) \leq F(n,c,m)$. Furthermore, the equality is attained if and only if $X$ is a variety of minimal degree.
\item[(2)] If $a_m (X) < F(n,c,m)$, then $ a_m (X) \leq G_{n+1} (n,c,m)$. Furthermore, the equality is attained if and only if $X$ is a del Pezzo variety.
\item[(3)] Suppose that $c \geq 3$. Then
\begin{equation*}
a_2 (X) = {{c+1} \choose {2}}-2~(=G_n (n,c,2))
\end{equation*}
if and only if either $d=c+2$ and ${\rm depth}(X)=n$ or else $d=c+3$ and ${\rm depth}(X)=n+1$.
\item[(4)] Suppose that $m \geq 3$. If $a_m (X) < G_{n+1} (n,c,m)$, then
\begin{equation*}
a_m (X) \leq G_n (n,c,m).
\end{equation*}
Furthermore, the equality is attained if and only if $d=c+2$ and ${\rm depth} (X)=n$.
\end{description}
\end{theorem}

\begin{proof}
(1) Proposition \ref{prop:degreebound} shows that $a_m (X) \leq F(n,c,m)$ (cf. (\ref{eq:4.9})) and the equality is attained only if $d=c+1$ (cf. (\ref{eq:4.5})). On the other hand, the equality holds if $d=c+1$ (cf. Definition and Remark \ref{lowdegreevarieties}.(A)).
\smallskip

\noindent (2) By (1) and Proposition \ref{prop:degreebound}, the condition $a_m (X) < F(n,c,m)$ implies that $d \geq c+2$ and hence
\begin{equation*}
a_m (X) \leq G_{n+1} (n,c,m)
\end{equation*}
(cf. (\ref{eq:4.9})). Also if $X$ is a del Pezzo variety, then the inequality turns into the equality (cf. Definition and Remark \ref{lowdegreevarieties}.(B)). Conversely, if $a_m (X) = G_{n+1} (n,c,m)$ then we have
\begin{equation*}
a_m (X)>H_3 (n,c,m)
\end{equation*}
by (\ref{eq:4.9}) and (\ref{eq:4.10}). Thus Proposition \ref{prop:degreebound} and (1) show that $d=c+2$. Now, (\ref{eq:4.5}) shows that $X$ is a del Pezzo variety.
\smallskip

\noindent (3) Suppose that
\begin{equation}\label{eq:4.15}
a_2 (X) ={{c+1} \choose {2}}-2 ~ (=H_3 (n,c,2) > H_4 (n,c,2)).
\end{equation}
Then $d = c+2$ or $c+3$ by (1) and Proposition \ref{prop:degreebound}. For $d=c+2$, we get ${\rm depth}(X)=n$ from (\ref{eq:4.3}). Now, we assume that $d=c+3$. Let $C \subset \P^{c+1}$ be a general curve section of $X$. Thus $C$ is a curve of degree $d=c+3$ such that
\begin{equation*}
a_2 (C) \geq a_2 (X) = {{c+1} \choose {2}}-2.
\end{equation*}
By (1) and (2), $a_2 (C)$ can be strictly larger than ${{c+1} \choose {2}}-2$ only if $d \leq c+2$. Thus, $a_2 (C)$ is equal to ${{c+1} \choose {2}}-2$ since $d=c+3$. Then Theorem \ref{thm:curvehomologicaltype} shows that $C$ is a linearly normal curve of arithmetic genus $2$. In particular, $C$ is an arithmetically Cohen-Macaulay curve (cf. Proposition \ref{prop:lowdegcurve}). This completes the proof that ${\rm depth}(X)=n+1$. Conversely, if either $d=c+2$ and ${\rm depth}(X)=n$ or else $d=c+3$ and ${\rm depth}(X)=n+1$, then $a_2 (X) ={{c+1} \choose {2}}-2$ by (\ref{eq:4.3}) and Definition and Remark \ref{lowdegreevarieties}.(D).
\smallskip

\noindent (4) By (1) and (2), the condition $a_m (X) < G_{n+1} (n,c,m)$ implies that either $d = c+2$ and ${\rm depth}(X) \leq n$ or else $d \geq c+3$. In the first case, we get $a_m (X) \leq G_n (n,c,m)$ with the equality if and only if ${\rm depth}(X)=n$ (cf. (\ref{eq:4.5})). In the second case, we get
\begin{equation*}
a_m (X) \leq H_3 (n,c,m) < G_n (n,c,m)
\end{equation*}
(cf. Proposition \ref{prop:degreebound} and (\ref{eq:4.12})). This completes the proof.
\end{proof}

\begin{remark}\label{rmk:analyizingtheproof}
From the proof of Theorem \ref{thm:firstandsecond}, one can try to solve Problem A in Section 1 by classifying all $X$ with $a_2 (X) = {{c+1} \choose {2}}+1-k$ and then extending the classification result to arbitrary $m \geq 3$. In this direction, Proposition \ref{prop:degreebound} shows that if $k \leq c$ and
\begin{equation*}
a_2 (X) = {{c+1} \choose {2}}+1-k ~ (=H_k (n,c,2) > H_{k+1} (n,c,2)),
\end{equation*}
then $d \leq c+k$. Thus one needs a structure theory for projective varieties of degree $c+k$. This was done in \cite{HSV} for $k=2$ (cf. Definition and Remark \ref{lowdegreevarieties}.(B))  but is still widely open even for $k=3$.
\end{remark}

\begin{remark}\label{rem:4.7}
Theorem \ref{thm:firstandsecond} shows that if $c \geq 2$ and $m \geq 3$, then
$$\delta_{n,c,m} (k) = \begin{cases} F(n,c,m) & \mbox{for $k=1$,}\\
                                     G_{n+1} (n,c,m) & \mbox{for $k=2$ and}\\
                                     G_n (n,c,m) & \mbox{for $k=3$.} \end{cases}$$
Also it can be shown from Definition and Remark \ref{lowdegreevarieties} that
\begin{equation*}
\delta_{n,c,m} (4) = \mbox{max} \{G_{n-1} (n,c,m), H_3 (n,c,m)\}.
\end{equation*}
More precisely,
\begin{equation*}
G_{n-1} (n,c,m) - H_3 (n,c,m) = \frac{(m+n-4)!}{n! (m-2)!} (m^2 +mn-n^2 -5m -n +6)
\end{equation*}
and hence
$$\delta_{n,c,m} (4) = \begin{cases} G_{n-1} (n,c,m) & \mbox{if $m^2 +mn-n^2 -5m -n +6 \geq 0$, and}\\
                                     H_3 (n,c,m)     & \mbox{otherwise.}\end{cases}$$
Moreover, if
\begin{equation*}
\delta_{n,c,m} (4) = G_{n-1} (n,c,m)
\end{equation*}
then $a_m (X) = \delta_{n,c,m} (4)$ if and only if $d=c+2$ and ${\rm depth}(X)=n-1$, and if
\begin{equation*}
\delta_{n,c,m} (4) = G_{n-1} (n,c,m)
\end{equation*}
then $a_m (X) = \delta_{n,c,m} (4)$ if and only if $d=c+3$ and ${\rm depth}(X)=n+1$.
\end{remark}
\vspace{0.5 cm}

\section{Projective invariants of quadratic embeddings}

\noindent The aim of this section is to investigate projective invariants of quadratic embedding of projective varieties having many quadratic equations.

In \cite{Z}, F. L. Zak defined several projective invariants of an embedded projective variety by using the higher secant varieties of its quadratic embedding. Also he established foundational works about them (eg. see Theorem \ref{thm:quadraticinvariants}). To be precise, let $X \subset \P^r$ be a nondegenerate projective variety of dimension $n$ and let
\begin{equation*}
Y := \nu_2 (X) \subset \P^N, \quad N= {{r+2} \choose {2}}-1,
\end{equation*}
be its quadratic embedding. Thus the span of $Y$ in $\P^N$ is a linear space of dimension $N-a_2 (X)$. For each $k \geq 1$, we denote by $S^k Y$ the variety swept out by the $k$-dimensional linear subspaces of $\P^N$ that are $(k+1)$-secant to $Y$. According to \cite{Z}, we consider the following projective invariants of $X$:\\

$s_k = {\rm dim}~S^k Y$

$k_2 = {\rm min}~ \{ k ~|~ S^k Y = \langle Y \rangle \}$

$\delta_k = s_{k-1}+n+1-s_k$ : the $k$-th secant deficiency of $Y$

$\ell_2 = {\rm max}~ \{ k ~|~ \delta_k = 0 \}$

$\delta^2  = \Sigma_{k = \ell_2 (X)+1} ^{k_2 (X)} ~ \delta_k$ : the total quadratic deficiency of $X$\\

Keeping the above notations in mind, we obtain the following

\begin{theorem}\label{thm:mainquadrics}
Let $X \subset \P^{n+c}$, $c>0$, be a nondegenerate projective variety of dimension $n$ and degree $d$.

\smallskip

\renewcommand{\descriptionlabel}[1]%
             {\hspace{\labelsep}\textrm{#1}}
\begin{description}
\setlength{\labelwidth}{13mm}
\setlength{\labelsep}{1.5mm}
\setlength{\itemindent}{0mm}

\item[(1)] The following conditions are equivalent:
\begin{enumerate}
\item[(i)] $X$ is a variety of minimal degree;
\item[(ii)] $a_2 (X)={{c+1} \choose {2}}$;
\item[(iii)] $\delta_{c+1} =1$.
\end{enumerate}
In this case, $\ell_2 = c$, $k_2 = n+c$ and $\delta_k = k-c$ for all $c+1 \leq k \leq c+n$.
\item[(2)] If $c \geq 2$ and $n \geq 2$, then the following conditions are equivalent:
\begin{enumerate}
\item[(i)] $X$ is a del Pezzo variety;
\item[(ii)] $a_2 (X)={{c+1} \choose {2}}-1$;
\item[(iii)] $\ell_2 = c+1$ and $\delta_{c+2} = 2$.
\end{enumerate}
In this case, $\ell_2 = c+1$, $k_2 = n+c$ and $\delta_k = k-c$ for all $c+2 \leq k \leq c+n$.
\item[(3)] If $c \geq 3$ and $n \geq 3$, then the following conditions are equivalent:
\begin{enumerate}
\item[(i)] $d=c+2$ and ${\rm depth}(X)=n$ or $d=c+3$ and ${\rm depth}(X)=n+1$;
\item[(ii)] $a_2 (X)={{c+1} \choose {2}}-2$;
\item[(iii)] $\ell_2 = c+1$, $\delta_{c+2} =1$ and $\delta_{c+3}=3$.
\end{enumerate}
In this case, $\ell_2 = c+1$, $k_2 = n+c$, $\delta_{c+2}=1$ and $\delta_k = k-c$ for all $c+3 \leq k \leq c+n$.
\end{description}
\end{theorem}
\smallskip

We begin with summarizing a part of main results in \cite{Z}:

\begin{theorem}[F. L. Zak, \cite{Z}]\label{thm:quadraticinvariants}
Let $X \subset \P^{n+c}$ be as in Theorem \ref{thm:mainquadrics}.

\smallskip

\renewcommand{\descriptionlabel}[1]%
             {\hspace{\labelsep}\textrm{#1}}
\begin{description}
\setlength{\labelwidth}{13mm}
\setlength{\labelsep}{1.5mm}
\setlength{\itemindent}{0mm}

\item[(1)] Let $\ell_2 \leq k < n+c$. Then $\delta_k < \delta_{k+1}$. In other words, $\delta_k$ is a strictly monotonous function in the interval $[\ell_2 , n+c]$.
\item[(2)] $\delta_k \begin{cases} =0 & \mbox{if $k \leq c$, and}\\
                            \leq k-c & \mbox{if $c+1 \leq k \leq c+n$.} \end{cases}$ \\
In particular, $\ell_2 \geq c$.
\item[(3)] If $\delta_k = k-c$ for some $c+1 \leq k \leq c+n$, then
\begin{equation*}
\delta_{k'} = k' -a
\end{equation*}
for all $k \leq k' \leq c+n$ and $k_2 = c+n$.
\item[(4)] \begin{equation*}
a_2 (X) = \delta^2 - (k_2 +1)(n+1)+ {{c+n+2} \choose {2}}.
\end{equation*}
\item[(5)] \begin{equation*}
a_2 (X) \quad \leq \sum_{k= \ell_2 +1} ^ {c+n} ~ \delta_k - {{n+1} \choose {2}}+{{c+1} \choose {2}}
\end{equation*}
\begin{equation*}
        \quad \quad \quad  \leq \delta^2 - {{n+1} \choose {2}}+{{c+1} \choose {2}},
\end{equation*}
where both inequalities turn into equalities if and only if $k_2 =c+n$.
\end{description}
\end{theorem}

\begin{proof}
See Theorem 3.1, 3.3, Corollary 3.4, Proposition 5.2, 5.3, Corollary 5.8 and Proposition 5.10 in \cite{Z}.
\end{proof}

\begin{corollary}\label{cor:inequality}
Let $X \subset \P^{n+c}$ be as in Theorem \ref{thm:mainquadrics} such that
\begin{equation*}
a_2 (X) = {{c+1} \choose {2}}-m \quad \mbox{for some $m \geq 0$.}
\end{equation*}

\smallskip

\renewcommand{\descriptionlabel}[1]%
             {\hspace{\labelsep}\textrm{#1}}
\begin{description}
\setlength{\labelwidth}{13mm}
\setlength{\labelsep}{1.5mm}
\setlength{\itemindent}{0mm}

\item[(1)] \begin{equation*}
0 \leq \sum_{k=c+1} ^{c+n} ~ (k-c-\delta_k) \leq m.
\end{equation*}
\item[(2)] If $m \leq n$, then $k_2 = n+c$ and $\delta^2 = {{n+1} \choose {2}} - m$. Furthermore, if $m < n$ then
\begin{equation*}
\delta_k = k-c  \quad \mbox{for all $c+m+1 \leq k \leq c+n$.}
\end{equation*}
\item[(3)] $m=0$ if and only if $\ell_2 =c$ (and hence $\delta_{c+1}=1$). In this case, $k_2 =n+c$ and
\begin{equation*}
\delta_k = k-c \quad \mbox{for all $c+1 \leq k \leq c+n$.}
\end{equation*}
\item[(4)] Suppose that $n \geq 2$. Then $m=1$ if and only if $\ell_2 = c+1$ and $\delta_{c+2} = 2$. In this case, $k_2 =n+c$ and
\begin{equation*}
\delta_k = k-c \quad \mbox{for all $c+2 \leq k \leq c+n$.}
\end{equation*}

\item[(5)] Suppose that $n \geq 3$. Then  $m=2$ if and only if $\ell_2 = c+1$, $\delta_{c+2} =1$ and $\delta_{c+3}=3$. In this case, $k_2 =n+c$ and
\begin{equation*}
\delta_k = k-c \quad \mbox{for all $c+3 \leq k \leq c+n$.}
\end{equation*}
\end{description}
\end{corollary}
\smallskip

\begin{proof}
(1) This follows immediately from Theorem \ref{thm:quadraticinvariants}.(2) and (5).
\smallskip

\noindent (2) Letting $u_k := k-c-\delta_k$, it holds by Theorem \ref{thm:quadraticinvariants}.(1) and (2) that
\begin{equation}\label{eq:5.1}
u_{c+1} \geq u_{c+2} \geq \cdots \geq u_{c+n} \geq 0.
\end{equation}
If $u_{c+n} \geq 1$, then (1) and (4.1) show that
\begin{equation*}
n \leq u_{c+1}+u_{c+2}+\cdots+ u_{c+n} \leq m.
\end{equation*}
Therefore we get $m=n$ and $u_{c+1}= u_{c+2}= \cdots= u_{c+n} =1$. In this case, the first inequality in Theorem \ref{thm:quadraticinvariants}.(5) turns into equality and so $k_2 = n+c$. If $u_{c+n}=0$, then we get $k_2 = n+c$ by Theorem \ref{thm:quadraticinvariants}.(3). Now, the equality $\delta^2 = {{n+1} \choose {2}} - m$ comes from Theorem \ref{thm:quadraticinvariants}.(5). Suppose that $m < n$. If $\delta_{c+m+1} \neq m+1$, then $\delta_{c+m+1} >0$ and so $u_{c+m+1} \geq 1$ by Theorem \ref{thm:quadraticinvariants}.(2). Therefore
\begin{equation*}
m+1 \leq u_{c+1}+u_{c+2}+\cdots+ u_{c+n}
\end{equation*}
by (\ref{eq:5.1}), which contradicts to (1). Consequently, it must be true that $\delta_{c+m+1} = m+1$. Then the last assertion follows by Theorem \ref{thm:quadraticinvariants}.(3).

\noindent (3) $(\Longrightarrow)$: By (2), we get $k_2 =n+c$ and $\delta_k = k-c$ for all $c+1 \leq k \leq c+n$. Therefore $\ell_2 =c$ by Theorem \ref{thm:quadraticinvariants}.(2).

$(\Longleftarrow)$: If $\ell_2 =c$, then $\delta_{c+1}>0$ and hence $\delta_{c+1}=1$ and $k_2 =n+c$ by Theorem \ref{thm:quadraticinvariants}.(2). Then by Theorem \ref{thm:quadraticinvariants}.(3), we get $\delta_k = k-c$ for all $c+1 \leq k \leq c+n$. Therefore $m=0$ by Theorem \ref{thm:quadraticinvariants}.(5).
\smallskip

\noindent (4) $(\Longrightarrow)$: By (2), we get $k_2 =n+c$, $\delta^2 = {{n+1} \choose {2}}-1$ and $\delta_k = k-c$ for all $c+2 \leq k \leq c+n$. Therefore $\delta_{c+1}=0$ and $\ell_2 =c+1$.

$(\Longleftarrow)$: If $\delta_{c+2} =2$, then $k_2 =n+c$ and $\delta_k = k-c$ for all $c+3 \leq k \leq c+n$ by Theorem \ref{thm:quadraticinvariants}.(3). Also $\delta_{c+1}=0$ since $\ell_2 =c+1$. Now, we get $m=1$ by Theorem \ref{thm:quadraticinvariants}.(5).

\noindent (5) $(\Longrightarrow)$: By (2), we get $k_2 =n+c$, $\delta^2 = {{n+1} \choose {2}}-2$ and $\delta_k = k-c$ for all $c+2 \leq k \leq c+n$. Also $\delta_{c+1}=0$ by (3). Therefore $\delta_{c+2}=1$.

$(\Longleftarrow)$: If $\delta_{c+3} =3$, then $k_2 =n+c$ and $\delta_k = k-c$ for all $c+3 \leq k \leq c+n$ by Theorem \ref{thm:quadraticinvariants}.(3). Also $\delta_{c+1}=0$ since $\ell_2 =c+1$. Now, we get $m=2$ by Theorem \ref{thm:quadraticinvariants}.(5).
\end{proof}
\smallskip

Now we give the \\

\noindent {\bf Proof of Theorem \ref{thm:mainquadrics}} For each of (1) $\sim$ (3), the equivalence of (i) and (ii) follows by Theorem \ref{thm:firstandsecond} and the equivalence of (ii) and (iii) is shown in Corollary \ref{cor:inequality}. \qed \\

\begin{example}\label{ex:VAMDdepthone}
Let $X \subset \P^{n+c}$, $c>0$, be an $n$-dimensional variety of almost minimal degree such that ${\rm depth}(X)=1$. Such an $X$ is always obtained as the image of an isomorphic linear projection of a smooth variety $\widetilde{X} \subset \P^{n+c+1}$ of minimal degree from a point. In particular, $X$ is smooth (cf. Theorem 1.1 in \cite{P}). Since $X$ is $3$-regular (cf. Theorem A in \cite{HSV}), it satisfies $2$-normality. This implies that the quadratic embedding of $X$ is that of $\widetilde{X}$. This observation and Theorem \ref{thm:mainquadrics}.(1) enable us to conclude that
$$\begin{cases}
\ell_2 (X) = \ell_2 (\widetilde{X})=c+2,\\
k_2 (X) = k_2 (\widetilde{X})= n+c+1, \\
\delta_k (X) = \delta_k (\widetilde{X}) = k-c-1 \quad \mbox{for all $c+2 \leq k \leq n+c+1$.}
\end{cases}$$
\end{example}
\smallskip

\begin{remark}
(1) Theorem \ref{thm:mainquadrics}.(1) and (2) are firstly shown at Corollary 5.8 and Proposition 5.10 in \cite{Z}, in which the statement of Proposition 5.10.(a.ii) in \cite{Z} should be replaced to that of our Theorem \ref{thm:mainquadrics}.(2.i).
\smallskip

\noindent (2) By Theorem \ref{thm:mainquadrics} and Example \ref{ex:VAMDdepthone}, we get the following Table 2 of the projective invariants of quadratic embeddings for some varieties of low degree.

\begin{table}[hbt]
\begin{center}
\begin{tabular}{|c|c|c|c|c|c|c|c|c|c|}\hline
${\rm deg}(X)$ & ${\rm depth}(X)$ & $\delta_{c+1}$ & $\delta_{c+2}$ & $\delta_{c+3}$ & $\cdots$ & $\delta_{c+i}$ & $\cdots$ & $\delta_{c+n}$ & $\delta_{c+n+1}$ \\\hline
$c+1$          & $n+1$            & $1$          & $2$              & $3$            & $\cdots$ & $i$            & $\cdots$ & $n$            & $0$              \\\hline
               & $1$              & $0$          & $1$              & $2$            & $\cdots$ & $i-1$          & $\cdots$ & $n-1$          & $n$              \\\cline{2-10}
$c+2$          & $n$              & $0$          & $1$              & $3$            &  $\cdots$ & $i$            & $\cdots$ & $n$            & $0$              \\\cline{2-10}
               & $n+1$            & $0$          & $2$              & $3$            & $\cdots$ & $i$            & $\cdots$ & $n$            & $0$              \\\hline
$c+3$          & $n+1$            & $0$          & $1$              & $3$            & $\cdots$ & $i$            & $\cdots$ & $n$            & $0$              \\\hline
\end{tabular}
\end{center}
\caption{Projective Invariants of Quadratic Embedding of $X \subset \P^{n+c}$}
\end{table}
\end{remark}

\begin{remark}
According to Theorem \ref{thm:mainquadrics}, one can ask the problem of classifying $n$-dimensional projective varieties
\begin{equation}\label{eq:5.2}
X \subset \P^{n+c} \quad \mbox{with $c \geq 4$ and $a_2 (X)={{c+1} \choose {2}}-3$}
\end{equation}
and to investigate their projective invariants of quadratic embeddings.
\smallskip

\noindent (1) By Corollary \ref{cor:quadraticinequality}, (\ref{eq:5.2}) implies that $d = c+2, c+3$ or $c+4$. Also by (\ref{eq:4.3}) and Corollary \ref{cor:boundarycases}, it follows that
\begin{enumerate}
\item[(i)] if $d=c+2$, then (\ref{eq:5.2}) holds if and only if ${\rm depth}(X)=n-1$
\end{enumerate}
and
\begin{enumerate}
\item[(ii)] if $d=c+4$, then (\ref{eq:5.2}) holds if and only if ${\rm depth}(X)=n+1$.
\end{enumerate}
Suppose that $d=c+3$ and let $C \subset \P^{c+1}$ be a general curve section of $X$. Then it can be shown that
\begin{equation*}
a_2 (C) ={{c+1} \choose {2}}-3
\end{equation*}
and $C$ is not linearly normal (see (\ref{eq:4.14}) and Proposition \ref{prop:lowdegcurve}). Therefore $C$ is one of the two curves in Table 1 for $k=3$. But it is not yet known which surfaces can take such curves as a hyperplane section.
\smallskip

\noindent (2) By a similar argument as in the proof of Corollary \ref{cor:inequality}, it can be shown that if $n \geq 3$ then $a_2 (X)={{c+1} \choose {2}}-3$ holds if and only if either 
\begin{enumerate}
\item[(i)] $\ell_2 = c+1$, $\delta_{c+2} =1$ and $\delta_{c+3}=2$ 
\end{enumerate}
or
\begin{enumerate}
\item[(ii)] $\ell_2 = c+2$ and $\delta_{c+3}=3$. 
\end{enumerate}
Thus it is an interesting question to ask if both of the above projective invariants of quadratic embeddings can occur.
\end{remark}
\vspace{0.5 cm}

\section{Hypersurfaces containing projective curves}

\noindent In this section, we study Problem A in Section 1 for $n=1$ and $m \geq c$. We will determine the value of $\delta_{1,c,m} (k)$ when $m \geq c$ and $1 \leq k \leq {{c} \choose {2}}+c$.

\begin{notation and remarks}\label{not-rem:5.1}
Recall that $\Xi_{1,c}$ is the set of all nondegenerate projective integral curves in $\P^{c+1}$, $A_{1,c,m} := \{ a_m (C) ~|~ C \in \Xi_{1,c} \}$ and for each $k \geq 1$, $\delta_{1,c,m} (k)$ is the $k$th largest member of $A_{1,c,m}$.
\smallskip

\noindent (A) Let $C \in \Xi_{1,c}$ be of degree $d$ and arithmetic genus $g$. Then for any $m \geq 2$,
\begin{equation*} 
a_m (C) = {{c+1+m} \choose {m}} - \{dm+1-g\} - h^1 (C,\mathcal{O}_C (m)) + h^1 (\P^{c+1},\mathcal{I}_C (m)).
\end{equation*}
Thus it is a hard question to say much about the value of $a_m (C)$ for arbitrary $m$. In this direction, the fact that ${\rm reg}(C) \leq d-c+1$ (see \cite{GLP}) is very important since it implies that
\begin{equation}\label{eq:a.2}
a_m (C) = {{c+1+m} \choose {m}} - \{dm+1-g\}
\end{equation}
for all $m \geq d-c$. In particular, (\ref{eq:a.2}) holds if $d \leq 2c$ and $m \geq c$ (see also Proposition \ref{prop:lowdegcurve} and Theorem \ref{thm:monotonousfunction}.(3)).
\smallskip

\noindent (B) For simplicity, we define the following integers:
\begin{equation*}
u_{c,g,d} (m) = {{c+1+m} \choose {m}}- \{ d m +1-g \}.
\end{equation*}
If $1 \leq k \leq {{c} \choose {2}}+c$, then there exist unique integers $g$ and $d$ such that
\begin{equation}\label{eq:a.3}
0 \leq g \leq c-1, \quad c+1+g \leq d \leq 2c \quad \mbox{and} \quad k={{d-c} \choose {2}}+(d-c-g).
\end{equation}
Moreover, the map
\begin{equation*}
\rho_m : \{ k \in \N ~|~ 1 \leq k \leq {{c} \choose {2}}+c \} \rightarrow  A_{1,c,m}, \quad k \mapsto u_{c,g,d} (m)
\end{equation*}
is injective and order-reversing for all $m \geq c$.
\end{notation and remarks}

We obtain the following

\begin{theorem}\label{thm:curvehigherdegree}
Let $m$ and $k$ be integers such that $m \geq c$ and $1 \leq k \leq {{c} \choose {2}}+c$. Then
\begin{equation*}
\delta_{1,c,m}(k) = {{c+1+m} \choose {m}}- \{ d m +1-g \}
\end{equation*}
where $d$ and $g$ are unique integers satisfying (\ref{eq:a.3}). Furthermore, if $m \geq c+1$ or if $m=c$ and $k < {{c} \choose {2}}+c$ then a curve $C \in \Xi_{1,c}$ satisfies $a_m (C) = \delta_{1,c,m} (k)$ if and only if it is of arithmetic genus $g$ and degree $d$.
\end{theorem}

\begin{proof}
Let $g$ and $d$ be integers such that $0 \leq g \leq c-1$ and $c+1+g \leq d \leq 2c$. Then a linearly normal smooth curve $\widetilde{C} \subset \P^{d-g}$ of genus $g$ and degree $d$ always exists. This implies that there is a nondegenerate projective smooth curve $C \subset \P^{c+1}$ of genus $g$ and degree $d$. Furthermore, ${\rm reg} (C) \leq c+1$ by Theorem \ref{thm:monotonousfunction}.(3) and hence
\begin{equation*}
a_m (C) = u_{c,g,d} (m)
\end{equation*}
for all $m \geq c$. That is, the set $M_{c,m}$ defined by
\begin{equation*}
M_{c,m} := \{ u_{c,g,d} (m) ~|~ 0 \leq g \leq c-1, ~ c+1+g \leq d \leq 2c \}
\end{equation*}
coincides with the subset
\begin{equation*}
\{ a_m (C) ~|~ C \in \Xi_{1,c},~ {\rm deg}(C) \leq 2c \}
\end{equation*}
of $A_{1,c,m}$. Also Notation and Remarks \ref{not-rem:5.1}.(B) says that $M_{c,m}$ consists of ${{c} \choose {2}}+c$ distinct integers and its $k$th largest member is $\rho_m (k)$. In particular, the smallest member of $M_{c,m}$ is
\begin{equation*}
\rho_m ({{c} \choose {2}}+c) = u_{c,0,2c} (m) = {{c+1+m} \choose {m}}-\{2cm+1\}.
\end{equation*}
On the other hand, Proposition \ref{prop:degreebound} says that if ${\rm deg}(C) \geq 2c+1$ then
\begin{equation*}
a_m (C) \leq H_{c+1} (1,c,m) =  {{c+1+m} \choose {m}}-\{(2c+1)m + 1 - c\}.
\end{equation*}
This completes the proof since $u_{c,0,2c} (m) - H_{c+1} (1,c,m) = m-c \geq 0$.
\end{proof}
\vspace{0.5 cm}

\bibliographystyle{plain}

\vskip 1 cm

\author{
 \begin{tabular}{ll}
Department of Mathematics \\
Korea University \\
Seoul 136-701 \\
Republic of Korea \\
{\it email: } euisungpark@korea.ac.kr
 \end{tabular}
}

\end{document}